\documentclass[12pt]{article}
\usepackage{amssymb}

\newtheorem{lemma}{Lemma}
\newtheorem{theorem}{Theorem}

\newtheorem{corollary}{Corollary}

\newcommand{\Z}{{\mathbb Z}}
\newcommand{\C}{{\mathbb C}}
\newcommand{\R}{{\mathbb R}}
\newcommand{\E}{{\mathbb E}}
\newcommand{\eps}{\varepsilon}

\newcommand{\B}{{\cal B}}

\title{Some properties of lower level-sets of convolutions}

\author{Ernie Croot \thanks{Supported by NSF grant DMS-1001111}}

\begin{document}

\maketitle

\abstract{In the present paper we prove a certain lemma about the structure of ``lower level-sets
of convolutions'', which are sets of the form $\{x \in \Z_N : 1_A*1_A(x) \leq \gamma N\}$ or of the form
$\{x \in \Z_N : 1_A*1_A(x) < \gamma N\}$, where $A$ is a subset
of $\Z_N$.  One result we prove using this lemma is that if $|A| = \theta N$ and $|A+A| \leq (1-\eps) N$, $0 < \eps < 1$, 
then this level-set contains an arithmetic progression of length at least $N^c$, $c = c(\theta, \eps,\gamma) > 0$.
It is perhaps possible to obtain such a result using Green's arithmetic regularity lemma (in combination with 
some ideas of Bourgain \cite{bourgain}); however, our method of proof allows us to obtain non-tower-type quantitative dependence
between the constant $c$ and the parameters $\theta$ and $\eps$.
For various reasons (discussed in the paper) one might think, wrongly, that such results would only 
be possible for level-sets involving triple and higher convolutions.} 
\bigskip


\section{Introduction}

There are many conditions that one can give on a subset $A \subseteq G$,
where $G$ is a finite additive abelian group,
guaranteeing that $A+A = \{a+b : a,b \in G\} = G$ or that $A+A$ is {\it nearly} all of 
$G$ (e.g. if $|A| > |G|/2$ or if the non-trivial Fourier
coefficients of the indicator function $1_A$ are all ``small'').  And one might wonder
whether there are some simple conditions on the set $A+A$ {\it itself} guaranteeing that
it is all of $G$, or at least a substantial proportion of $G$; for example, is
there a particular small set $S$  such that if 
we know that $A+A$ contains $S$ and $|A| > |G|/\log |G|$, say, 
then $A+A$ must  essentially be all of $G$?  
In the present paper we develop some related results.  The key idea behind most of
them is a lemma (actually, Corollary \ref{main_lemma}) on the structure of level-sets of 
convolutions given in a later section.  It may be possible to obtain versions of some of
our results (particularly Theorem \ref{longaps} and its corollaries) 
using Green's arithmetic regularity lemma \cite{green_regular} in combination
with ideas of Bourgain on Bohr neighborhoods \cite{bourgain}; however, such techniques
typically give weaker quantitative bounds than what we produce in the present paper.

In order to discuss some of these results, we will need some notation:  
suppose that $G$ is a finite group and that $g : G  \to \C$.
We define the expectation operator
$$
\E g\ =\ \E_x g\ :=\ |G|^{-1} \sum_{x \in G} g(x);
$$
we define for an additive abelian group $G$ the (unnormalized) convolution $f*g$ of 
two functions $f,g : G \to \C$ to be
$$
f*g(x)\ :=\ \sum_{a+b = x \atop a,b \in G} f(a)g(b)\ =\ |G|\cdot \E_{a \in G} f(a) g(x-a);
$$
given $f : G \to \C$ we define the (unnormalized) Fourier transform at $\chi \in \hat G$ to be
$$
\hat f(\chi)\ :=\ \sum_{x \in G} f(x) \chi(x)\ =\ |G| \cdot \E_x f(x) \chi(x);
$$
and lastly, we say that a function $f : G \to \C$ is $\alpha$-uniform if 
$$
\max_{\chi \in \hat G} |\hat f(\chi)|\ \leq\ \alpha |G|.
$$
An easy consequence of the triangle inequality and the linearity of the Fourier
transform is that if $f_1,...,f_k$ are
$\alpha_1,...,\alpha_k$-uniform, respectively, then their sum $f_1+\cdots + f_k$
is $\alpha_1 + \cdots + \alpha_k$-uniform.

Our first result is along the lines of what we described above, except that we replace the
condition that $A+A$ contains $S$ with the condition that $1_A*1_A(s)$ is ``large'' for all
$s \in S$ -- such a condition is often easier and more natural to work with than having
$A+A$ contain $S$:

\begin{theorem} \label{A_convolve_theorem}  Suppose $G$ is an additive abelian group with
$|G| = N$, and suppose that $0 < \theta \leq 1$ and $\delta, \eps > 0$.
Then, there exists a set $S \subseteq G$ satisfying 
$$
|S|\ \ll\ \eps^{-2} \delta^{-6} \theta^{-10} (\log N - \log(\delta \theta \eps)),
$$ 
such that if $A \subseteq G$, $|A| =\theta_0 N \geq \theta N$ 
satisfies $1_A*1_A(x) > \delta \theta_0^2 N$ for every $x \in S$, then
$|A+A| > (1-\eps)N$.
\end{theorem}

\noindent {\bf Note:}  The expected value of $1_A*1_A$ is $\theta_0^2 N$; so, the condition $1_A*1_A(x) > \delta \theta_0^2 N$
is just requiring that the convolution be more than $\delta$ times as big as this expected value.
\bigskip

This theorem is not far from best-possible in that $|S|$ needs to be $\Omega(\log N)$ in
order for the conclusion of the theorem to hold, as the following result demonstrates for 
$\theta \geq 1/3$:

\begin{theorem} \label{A_sum_limits}  For every sufficiently large prime $N$,
and every set $S \subseteq \Z_N$ of
size at most $(\log N)/2$, there exists a set $A \subseteq \Z_N$ of size $|A| \geq N/3$
such that $|A+A| < 2N/3$ and such that $1_A*1_A(x) > N/6$ for every $x \in S$.
\end{theorem}
\bigskip

Not only is it possible to show that {\it there exists} a set $S$ having the requisite properties
given in Theorem \ref{A_convolve_theorem}, but, in fact, if we allow $|S|$ to be somewhat larger 
than Theorem \ref{A_convolve_theorem} requires then we can show that {\it any} set $S$ whose 
non-trivial Fourier coefficients are sufficiently ``small'' 
(indicated below) will do.  

\begin{theorem} \label{A_sum_pseudo}  Suppose that $G$ is an additive abelian group
satisfying $|G| = N$, and suppose that $\delta, \eps, \theta > 0$ are parameters
that we allow to depend on $N$.  Let $S \subseteq G$ be a set such that 
$$
\max_{\chi \in \hat G \atop \chi \neq \chi_0} |\hat 1_S(\chi)|\ <\  
(\delta^3 \theta^{4.5} \sqrt{\eps}/512 \pi)^{2^{16}\delta^{-6}\theta^{-10}\epsilon^{-1}+1}|S|,
$$  
where $\chi_0$ denotes the principal character.
Then, if $A \subseteq G$ satisfies $|A| \geq \theta N$
and $1_A*1_A(x) > \delta \theta^2 N$
for every $x \in S$,   we will have that $|A+A| \geq (1-\eps) N$.
\end{theorem}

See \cite{ajtai} for some explicit constructions of small sets $S$, all of whose non-trivial Fourier
coefficients are small.
\bigskip

The last set of results we prove are rather different from the ones listed 
above and pertain to the existence of long arithmetic progressions and other structures 
in level-sets of the convolution $1_A*1_A(x)$.  By {\it level-set} here we mean a set having the form
$\{ x \in G\ :\ 1_A*1_A(x) \leq \gamma |G|\}$ or having the form
$\{x \in G\ :\ 1_A*1_A(x) < \gamma |G|\}$.   As these results make use of Bohr neighborhoods,
now is a good time to define them:
\bigskip

\noindent {\bf Definition.}  Suppose that $\Lambda := \{\chi_1,...,\chi_k\} \subseteq \hat G$ and that
$\eps > 0$.  Then, the Bohr neighborhood of radius $\eps$ determined by $\Lambda$ is defined
to be the set
$$
\B(\Lambda, \eps)\ :=\ \{ x \in G\ :\ {\rm for\ }\ i=1,...,k,\ | 1 - \chi_i(x)| \leq \eps \}. 
$$
The {\it dimension} of a Bohr neighborhood is the least number of places $\chi_i$ 
needed to define the set; so, this $\B(\Lambda,\eps)$ we wrote down 
has dimension at most $k$.
\bigskip

Our first result along these lines is stated
as follows:

\begin{theorem} \label{longaps} Suppose that $G$ is an additive abelian group satisfying
$|G| = N$; suppose that $A \subseteq G$, $|A| = \theta N$; and suppose that for $\delta > 0$,
\begin{equation} \label{upper_condition}
|\{x \in G\ :\ 1_A*1_A(x) < \delta^3 \theta^6 N/128 \}|\ \geq\ \eps N.
\end{equation}
Then, we have that the level-set
\begin{equation} \label{lowerset}
\{ x \in G\ :\ 1_A*1_A(x) < \delta \theta^2 N\}
\end{equation}
contains a translate of a Bohr neighborhood of dimension at most 
$2^{16}\delta^{-6}\theta^{-10}\epsilon^{-1}+1$ 
and radius $\delta^3 \theta^{4.5} \sqrt{\eps}/128$.  

Furthermore, if $G = \Z_N$, where $N$
is prime, then using the fact that large Bohr neighborhoods always contain long arithmetic progressions, we have 
in this case that (\ref{lowerset}) contains an arithmetic progression of length at least $N^c$, where 
$c= c(\theta, \eps, \delta)$.
\end{theorem}

It is relatively straightforward to use standard Fourier arguments to deduce that 
the level-set (\ref{lowerset}) contains {\it most} of a translate of a Bohr neighborhood;
and, using ideas due originally to Bogolyubov \cite{bogolyubov} one can deduce that
the triple-convolution analogue of (\ref{lowerset}) -- e.g. $\{x \in G : 1_A*1_A*1_A(x)
< \delta \theta^3 N^2\}$ -- contains a complete shifted Bohr neighborhood
(see also \cite{green} and \cite{sanders}).  Using Green's arithmetic regularity 
lemma \cite{green_regular} in combination with some techniques in the manipulation
of Bohr neighborhoods due to Bourgain \cite{bourgain}, it is perhaps possible to prove a version of
Theorem \ref{longaps}, but with worse bounds on the radius and 
dimension of the Bohr neighborhood (due to the fact tower-type dependencies of certain parameters
inherent in the regularity lemma).

By replacing the condition (\ref{upper_condition}) with simply an upper bound for $|A+A|$
we arrive at the following corollary:

\begin{corollary}  \label{long_aps_corollary} Suppose that $A \subseteq G$, $|A| = \theta N$, and
$|A+A| \leq (1-\eps )N$.  Then, we deduce that the set (\ref{lowerset})  contains the same translate of the
Bohr neighborhood indicated by Theorem \ref{longaps} (and the associated arithmetic progression
of length at least $N^c$ in the case $G = \Z_N$, $N$ prime).
\end{corollary}

Interestingly, a consequence of a construction of Ruzsa \cite{ruzsa} is that there exist
sets $A$ for which the complement of this level-set cannot contain progressions
of length larger than $\exp( (\log N)^{2/3+\eps})$.  So, there is something seemingly 
paradoxical going on:  the ``lower level sets'' $\{ x \in \Z_N : 1_A*1_A(x) \leq \delta \theta^2 N\}$
{\it always} contain power-of-$N$-length arithmetic progressions provided $|A+A|$ isn't too large,
whilst the ``upper level sets'' $\{ x \in \Z_N : 1_A*1_A(x) \geq \delta \theta^2 N\}$ {\it sometimes}
do not.

From (the contrapositive of) Corollary \ref{long_aps_corollary} one can immediately deduce the following 
additional corollary, which gives another way to think about the result:

\begin{corollary}  Suppose that $A \subseteq \Z_N$ satisfies $|A| = \theta N$, where $N$ is prime.  
Furthermore, suppose that for every $x \in \Z_N$  satisfying
$1_A*1_A(x) \geq \delta \theta^2 N$ and for every $d \in \Z_N$, $d \neq 0$,
we have that there exists $0 < t \leq N^c$, $c = c(\theta, \eps, \delta) > 0$, such that 
$1_A*1_A(x + dt) \geq \delta \theta^2 N$.  
(In other words, the gaps along APs of common difference $d=1,...,N-1$ 
where $1_A*1_A(x)$ is ``large'', are all bounded from above by $N^c$.)  Then, 
$|A+A| \geq\ (1-\eps) N$.
\end{corollary}
\bigskip

The remainder of the paper is organized as follows:  in section \ref{future_section} 
we list some conjectures that, if true, would give much stronger structural information 
about lower level-sets than Theorem \ref{longaps} provides; in section \ref{auxiliary_section}
we list out two technical lemmas used throughout the rest of the paper; in section 
\ref{key_section} we state and prove the main lemma of the paper, and discuss some immediate 
consequences of it; in section \ref{proofs_section} we prove Theorems \ref{A_convolve_theorem}, 
\ref{A_sum_limits}, \ref{A_sum_pseudo}, \ref{longaps};  finally, the remaining sections are 
devoted to acknowledgements and the bibliography.

\section{Future directions} \label{future_section}

In this section we list some conjectures motivated by Theorem \ref{longaps} that would be natural
``next steps'' for where to continue this work.

We begin by noting that Theorem \ref{longaps} shows that certain level-sets always contain
long arithmetic progressions (when the group is $G = \Z_N$, $N$ prime), and one might
wonder wether {\it complements} of sumsets always contain such long arithmetic progressions (which
is itself the level-set $\{x : 1_A*1_A(x) = 0\}$); in other words, in order to get the long-progressions
conclusion in Theorem \ref{longaps}, must one necessarily work with level-sets of the form 
$\{x : 1_A*1_A(x) < \delta \theta^2 N\}$,
where $\delta > 0$?  This motivates the following conjecture:
\bigskip

\noindent {\bf Conjecture 1.}  For every sufficiently large prime $N$ there exists a set 
$A \subseteq \Z_N$, $|A| > N/4$, say, such that 
$|A+A| \leq 99 N/100$ (say), and such that the longest arithmetic progression in the complement of 
$A+A$ has length at most $N^{o(1)}$.
\bigskip

Perhaps something like Ruzsa's construction \cite{ruzsa} can be made to prove this conjecture; 
however, the author has not yet been able to do so.
\bigskip

The author also cannot immediately see any reason why the results in Theorem \ref{longaps}
could not be made much stronger; perhaps, in fact, the following is true:
\bigskip

\noindent {\bf Conjecture 2.}  The following holds for certain absolute constants
$0< c_1,c_2,c_3,c_4 < 1$ and primes $N$ sufficiently large:  suppose that $A \subseteq \Z_N$,  $|A| > N \exp(-(\log N)^{c_1})$
and $|A+A| \leq N(1 - \exp(-(\log N)^{c_2}))$.  Then, the lower level-set $\{x \in \Z_N : 1_A*1_A(x) < N \exp(-(\log N)^{c_3})\}$ 
contains an arithmetic
progression of length at least $\exp( (\log N)^{c_4})$.
\bigskip

An even stronger conjecture would be that, up to some small error, the
level-set above is, in fact, {\it covered} by disjoint translates of some given large Bohr neighborhood.
One formulation of such a conjecture is as follows:
\bigskip

\noindent {\bf Conjecture 3.}  The following holds for certain absolute constants $0 < c_1,c_2,c_3,c_4, c_5 < 1$ and
primes $N$ sufficiently large:  suppose that $A \subseteq G$, where $G$ is a finite abelian group
of size $N$, where  $|A| > N \exp(-(\log N)^{c_1})$ and 
$|A+A| \leq N (1 - \exp(-(\log N)^{c_2}))$.  Then, there exists a Bohr neighborhood $\B$ of dimension
at most $(\log N)^{c_3}$ and size $|\B| \geq N \exp(-(\log N)^{c_4})$, and a set of translates
$t_1,...,t_k$ with $\B + t_i$ all disjoint, such that if we let 
$$
S\ :=\ \cup_{i=1}^k (t_i + \B),
$$
then the symmetric difference between $S$ and the lower level-set
$\{x \in \Z_N : 1_A*1_A(x) < N \exp(-(\log N)^{c_5})\}$ has size at most $|S|/100$.
\bigskip

The sort of approach that {\it might} work here would be to combine ideas from the 
present paper with those from \cite{croot} and \cite{sanders}; however, the
author currently cannot see how to do this.

Besides being interesting problems in their own right, the two last conjectures above 
{\it could} perhaps be used to deduce good upper bounds on the largest subset $A \subseteq \Z_N$
having no non-trivial solutions $x,y,z$ to a given 
linear equations $a_1 x + a_2 y + a_3 z \equiv 0 \pmod{N}$, where $a_1 + a_2 + a_3 \equiv 0 \pmod{N}$,
thereby sharpening the already remarkable results of Sanders \cite{sanders2},
to achieve a similar success for such problems as was recently done by Schoen and Shkredov 
in their sensational paper \cite{schoen} 
for equations involving six or more variables (that produced upper bounds for $|A|$
of the general form of the Behrend \cite{behrend} \cite{elkin} \cite{wolf}
bound for particular choices of the $a_i$'s) and by Bloom \cite{bloom} for equations involving 
four variables and higher (in this paper he beautifully generalized Sanders's proof \cite{sanders2}).  
Both of these papers \cite{bloom} and \cite{schoen} make use of ideas from
\cite{croot2}, \cite{sanders} and \cite{sanders2}.

There might also be a way to use a proof of Conjecture 3 to improve upon the breakthrough results of 
Bateman and Katz \cite{bateman}.

Here is a rough idea of exactly {\it how} the last two conjectures above might be applicable
to problems about solutions to linear equations:  it is easiest to describe this in the case where
$G = \Z_3^n$, so let us make this assumption; and, let us suppose that $A \subseteq G$ has
no solutions to $x + y - 2z = 0$ -- that is, no solutions to $x + y + z = 0$, since $-2 = 1$ in $\Z_3$.
It follows that $-A$ is a subset of the level-set $\{x \in G : 1_A*1_A(x) \leq 1\}$.  If we knew
that this level-set were approximately the disjoint union of translates of some large Bohr neighborhood
$\B$, then one of those translates $t+\B$ should intersect $-A$ in many elements; indeed, one
would expect that $(-A) \cap (t+\B)$ has a higher density in $t+\B$ than $-A$ does in $G$.
This is exactly what we need in order to implement a ``density increment
strategy'' for showing that $|A|$ must be rather small, as was first done by Roth \cite{roth}.

\section{Some auxiliary lemmas} \label{auxiliary_section}

In this section we list some basic lemmas (and prove two of them) that we use in later sections
in the proofs of our main theorems.

\begin{lemma} \label{convolutions_uniformity_lemma}  Suppose that $G$
is an additive abelian group with $|G| = N$, and that 
$g, h : G \to [0,1]$, $\E g = \theta$.  Then, if $g - h$ is $\delta_1$-uniform we will have
$$
\sum_{x \in G} |g*g(x) - h*h(x)|^2\ \leq\ \delta_1^2 (4\theta + \delta_1) N^3.
$$
\end{lemma}

\noindent {\bf Proof of the Lemma.}  First note that since
$g - h$ is $\delta_1$-uniform and since $\E(g) = \theta$, we have that 
$\E(h) \leq \theta + \delta_1$; and so, since $g,h : \Z_N \to [0,1]$ this
implies that
\begin{equation} \label{me1}
\sum_{x \in G} (g+h)(x)^2\ \leq\ (4 \theta + \delta_1) N.
\end{equation}

And now from Parseval's identity we then have:
\begin{eqnarray*}
\sum_{x \in G} |g*g(x) - h*h(x)|^2\ &=&\ \sum_{x \in G} |(g-h)*(g+h)(x)|^2  \\
&=&\  N^{-1} \sum_{\chi \in \hat G} |\widehat{(g - h)}(\chi)|^2 |\widehat{(g+h)}(\chi)|^2 \\
&\leq&\ \delta_1^2 N \sum_{\chi \in \hat G} |\widehat{(g+h)}(\chi)|^2 \\
&=&\ \delta_1^2 N^2 \sum_{x \in G} (g+h)(x)^2,
\end{eqnarray*}
which, in combination with (\ref{me1}), proves the lemma.
\hfill $\blacksquare$
\bigskip

\begin{lemma} \label{Fourier_upper_bound}
Suppose that $h : G \to [0,1]$, where $G$ is an additive abelian group
satisfying $|G| = N$.  Then, if we place
the Fourier coefficients in order from largest to smallest, 
\begin{equation} \label{respect_me}
|\hat h(\chi_1)|\ \geq\ |\hat h(\chi_2)|\ \geq\ \cdots\ \geq\  |\hat h(\chi_N)|,
\end{equation}
where $\{\chi_1,...,\chi_N\} = \hat G$, we will have that
$$
|\hat h(\chi_k)|\ \leq\ N \sqrt{\E(h)/k}.
$$
\end{lemma}

\noindent {\bf Note:}   Because we may have $|\hat h(\chi_i)| = |\hat h(\chi_{i+1})|$ for
some $i=1,...,N$, the order of the $\chi_i$'s is not necessarily well-defined.  
However, the conclusion of the lemma does not depend on this choice; furthermore,
for the applications of this lemma in later sections the choice of ordering of the
$\chi_i$'s does not matter so long as they respect (\ref{respect_me}).
\bigskip

\noindent {\bf Proof of the lemma.}  We have from Parseval's identity that
$$
k |\hat h(\chi_k)|^2\ \leq\ \sum_{\chi \in \hat G} |\hat h(\chi)|^2\ \leq\ \E(h) N^2.
$$
Solving for $|\hat h(\chi_k)|$, the lemma follows.
\hfill $\blacksquare$
\bigskip

Finally, we will need the following lemma, which can be found in \cite[sec. 4.4]{taovu}, 
that gives a lower bound on the cardinalities of certain Bohr neighborhoods:

\begin{lemma}  \label{bohr_lemma} Suppose $G$ is an additive abelian group satisfying $|G| = N$.
If $\Lambda \subseteq \hat G$, $|\Lambda| = d$, then 
for $r \in [0,2]$ we have $|\B(\Lambda, r)| \geq (r/2\pi)^d N$.  Furthermore, 
if $G = \Z_N$ where $N$ is prime then this Bohr neighborhood contains an arithmetic progression 
of size at least $r N^{1/d}/2\pi$. 
\end{lemma}

\section{The key lemma and its proof} \label{key_section}

Before we can state the main lemma, we need to define the notion of 
a ``generalized convolution'' for finite groups:  suppose that $G$ is a finite
group (possibly non-abelian) where the operation is written multiplicatively,
and suppose that $T : G \times G \to G$.  
Then, for two functions $f,g : G \to \C$ we define the {\it $T$-convolution of 
$f$ with $g$} to be 
$$
f*_Tg(x)\ =\ \sum_{a,b \in G \atop T(a,b) = x} f(a) g(b).
$$  
Associated with this generalized convolution, we define a parameter
$$
\kappa = \kappa(T)\ :=\ \max(\kappa_1, \kappa_2),
$$
where
\begin{eqnarray*}
\kappa_1\ &:=&\ \max_{x \in G} \max_{a \in G} |\{ b\in G\ :\ T(a,b) = x\}|;\ {\rm and\ }\\
\kappa_2\ &:=&\ \max_{x \in G} \max_{a \in G} |\{ b \in G\ :\ T(b,a) = x\}|.
\end{eqnarray*}

Note that in the case where $T(a,b) = ab$, the convolution $f*_T g(x)$ coincides with
the usual group convolution $f*g(x) = \sum_{a,b \in G \atop ab = x} f(a)g(b)$; 
furthermore, the parameter $\kappa = \kappa_1 = \kappa_2 = 1$ in this case.

Our main lemma is given as follows:

\begin{lemma} \label{main_lemma0}  Suppose that $\| \cdot \|$ is a norm 
on the space of functions $f : G \to \C$, where $G$ is a finite group (possibly
non-abelian) of size $N$, and that $T : G \times G \to G$ has associated parameter $\kappa$ as
defined above.
Further, suppose that $\delta_1, \delta_2 > 0$ are parameters we allow to depend on $N$,
and that $A \subseteq G$ satisfies $|A| = \theta N \geq \delta_2 N > 0$.  Then, 
there exists a function $f : G \to [0,1]$ satisfying $\| f - 1_A\| \leq \delta_1$ such that
for every $B, C \subseteq G$ with $\| 1_B - 1_A\|, \| 1_C - 1_A\| \leq \delta_1$, we
have that
$$
1_B*_T1_C(x)\ \leq\ \delta_2^{-2} f*_Tf(x) + 2\kappa \delta_2 N.
$$
\end{lemma}

In order to make much use of this lemma, it seems that the choice of norm should
somehow be related to the convolutions $1_B*_T 1_B$ and $f*_T f$.  In the case where
$G$ is an additive abelian group and $T(a,b) = a+b$, Lemma \ref{convolutions_uniformity_lemma} implies that 
there is a natural choice for the norm having
this property, namely we can use $\| g \| = N^{-1} \max_{\chi \in \hat G} |\hat g(\chi) |$, $N = |G|$.
This now brings us to the following immediate corollary of Lemma \ref{main_lemma0}:

\begin{corollary} \label{main_lemma}  Suppose that $G$ is an additive abelian group 
with $|G| = N$.  Suppose that $\delta_1, \delta_2 > 0$ are parameters we allow to depend on 
$N$, and that $A \subseteq G$ 
satisfies $|A| = \theta N \geq \delta_2 N > 0$.
Then, there exists a function $f : G \to [0,1]$
such that $f - 1_A$ is $\delta_1$-uniform and such that for every
$B, C \subseteq G$ having the properties that both $1_B - 1_A$ and $1_C - 1_A$ are 
$\delta_1$-uniform, we have
$$
1_B*1_C(x)\ \leq\ \delta_2^{-2} f*f(x) + 2 \delta_2 N.
$$
\end{corollary}

We now will try to give some idea of what Lemma \ref{main_lemma0} is saying
by considering the special case of the above corollary with $G = \Z_N$:  fix a subset $A \subseteq \Z_N$ of size $\theta N$
and then consider all the other sets $B \subseteq \Z_N$ such that $1_B - 1_A$ is 
``highly uniform'' -- that is, $1_B - 1_A$ is $\delta_1$-uniform, where $\delta_1 > 0$ is 
``small''.  Lemma \ref{convolutions_uniformity_lemma} then implies that the 
convolutions $1_B*1_B$ and $1_A*1_A$ are ``close'' to one another in an ${\rm L}^2$ sense; but they need not be close
in an ${\rm L}^\infty$ sense, and in fact they cannot be in general.  A good example to demonstrate
the point is to consider the case where $N \equiv 3 \pmod{4}$ is a prime number, and where
$A = \{x^2 \pmod{N}\ :\ 1 \leq x \leq N-1\}$, which has size about $N/2$.
This set has the property that $1_A*1_A(0) = 0$, while $1_A*1_A(x) \sim N/4$ for $x \neq 0$.
Furthermore, all the non-zero Fourier coefficients of $1_A$ are ``small'';  indeed,
for $a \not \equiv 0 \pmod{N}$ we have from Gauss sum estimates that 
$|\hat 1_A(a)| \ll \sqrt{N}$.  If we now define the set $A_t := A + t$ then we likewise will have that
$1_{A_t}*1_{A_t}(2t) = 0$, and that $|\hat 1_{A_t}(a)| \ll \sqrt{N}$ for $a \not \equiv 0 \pmod{N}$.
In particular, the level-sets $\{x \in \Z_N : 1_A*1_A(x) < N/8\}$ and 
$\{x \in \Z_N : 1_{A_t}*1_{A_t}(x) < N/8\}$, are disjoint for $t \not \equiv 0 \pmod{N}$
(the first level-set here is $\{0\}$, while the
second is $\{2t\}$), even though $1_A - 1_{A_t}$ is $O(1/\sqrt{N})$-uniform; in addition,
$$
\| 1_A*1_A - 1_{A_t}*1_{A_t}\|_\infty\ =\ \max_{x \in \Z_N} |1_A*1_A(x) - 1_{A_t}*1_{A_t}(x)|\ \sim\ N/4,
$$
which is rather large.  

It would seem that this is pretty much all that one can say on the intersection of 
level-sets; but Corollary \ref{main_lemma} says that if
the lower level-sets we are considering are all rather large (instead of just a single element
as in the example involving the squares mod $N$), then in fact they all have large intersection.  The following
corollary of the Corollary \ref{main_lemma} gives a quantitative version of this fact:

\begin{corollary}   Suppose $G$ is an additive abelian group satisfying $|G| = N$, and 
fix a subset $A \subseteq G$, $|A| = \theta N$.  Then, letting $0 < \gamma \leq 1$ and letting
$$
I\ :=\ \bigcap_{B \subseteq G \atop 1_A - 1_B\ {\rm is\ }\delta-{\rm uniform}} \{ x \in G\ :\ 1_B*1_B(x) \leq \gamma \theta^2 N\},
$$
we have that
$$
|I|\ \geq\ |\{ x \in G\ :\ 1_A*1_A(x) \leq \gamma^3 \theta^6 N/128 \}| - 2^{14} \gamma^{-6} \theta^{-12} \delta^2 (4\theta + \delta) N.
$$

If we furthermore suppose that $|A+A| \leq (1-\eps) N$ then we immediately deduce from this corollary that
$$
|I|\ \geq\ (\eps - 2^{14} \gamma^{-6} \theta^{-12} \delta^2 (4\theta + \delta)) N.
$$
For fixed $\theta, \eps, \gamma > 0$, then, we see that if $\delta > 0$ is sufficiently small in terms of 
$\theta, \eps$ and $\gamma$, we must have that $|I|\ \gg\ \eps N$.
\end{corollary}
\bigskip

\noindent {\bf Proof of the Corollary.}  Let $\delta_1 = \delta$, $\delta_2 = \gamma \theta^2/4$, and then let $f : G \to [0,1]$ be the 
function given by our Corollary \ref{main_lemma}.  We have then for every $x \in G$ such that 
$f*f(x) \leq \gamma^3 \theta^6 N/64$ and for every set $B \subseteq G$ where $1_A - 1_B$ is $\delta$-uniform,
$$
1_B*1_B(x)\ \leq\ \delta_2^{-2} f*f(x) + 2 \delta_2 N\ \leq\ 3\gamma \theta^2 N/4.
$$

Next we apply Lemma \ref{convolutions_uniformity_lemma} using $g = 1_A$, $h = f$, 
and deduce that if we let $S$ be the
set of all $x \in G$ such that $1_A*1_A(x) \leq \gamma^3 \theta^6 N/128$, and $T \subseteq S$
be those $x\in S$ where $f*f(x) > \gamma^3 \theta^6 N/64$, then
$$
|T| (\gamma^3 \theta^6 N/128)^2\ \leq\ \delta_1^2 (4\theta + \delta_1) N^3;
$$
that is,
$$
|T|\ \leq\ 2^{14} \gamma^{-6} \theta^{-12} \delta^2 (4 \theta + \delta) N,
$$
which completes the proof of the Corollary. \hfill $\blacksquare$

\subsection{Proof of Lemma \ref{main_lemma0}} 

The proof of this lemma iterates on single places $x \in G$ where
some convolution $1_B*_T 1_C(x)$ is ``large'', which makes it somewhat like
the Dyson $e$-transform and also the Katz-Koester Lemma \cite{katz}.

We begin by constructing a sequence of functions
$f_1,f_2,...$ according to the following algorithm:

\begin{enumerate}
\item Set $f_1 := 1_A$, and set $j := 1$.

\item Suppose we have constructed $f_j$.  If for every pair of
sets $B,C \subseteq G$
such that $\| 1_A-1_B\|, \|1_A - 1_C\| \leq \delta_1$ we have that
$$
1_B*_T 1_C(x)\ \leq\ f_j*_T f_j(x) + 2\kappa \delta_2 N\ {\rm for\ every\ }x \in G,
$$
then we STOP.

\item Otherwise, there exist sets $B,C\subseteq G$ for which 
$\|1_A-1_B\|, \|1_A-1_C\| \leq \delta_1$, and for which there exists $x \in \Z_N$
satisfying 
$$
1_B*_T 1_C(x)\ >\ f_j*_T f_j(x) + 2\kappa \delta_2 N.
$$
Given such sets $B,C$ we either set $f_{j+1} := f_j + 1_B$ or $f_{j+1} := f_j + 1_C$,
according to which of these two possibilities makes $f_{j+1}$ have the larger
support (if there is a tie in the size of the support, simply choose $f_{j+1} := f_j + 1_B$).
Then, set $j \leftarrow j+1$.

\item And then we loop back to the second step.
\end{enumerate}

Let us see that this procedure must eventually terminate:  first, we note that
$1 \leq f_j(x) \leq j$ for all $x \in {\rm supp}(f_j)$ and for all $j \geq 1$.  Given $f_j$, if
there exist sets $B,C \subseteq G$ as in the third step, we must have that there exists 
$x \in G$ satisfying
$$
1_B*_T 1_C(x)\ >\ f_j*_T f_j(x) + 2\kappa \delta_2 N.
$$
Thus, there are more than $f_j*_T f_j(x) + 2\kappa \delta_2 N$ pairs
$$
(b,c)\ \in\ B \times C\ {\rm with\ } T(b,c)=x,
$$ 
while there are at most $f_j*_T f_j(x)$ pairs 
$$
(b,c)\ \in\ {\rm supp}(f_j) \times {\rm supp}(f_j)\ {\rm with\ } T(b,c)=x.
$$ 
It follows that there are more than $2\kappa \delta_2 N$ pairs $(b,c) \in B \times C$ 
with $T(b,c) =x$ for which either $b$ or $c$ fails to belong to 
${\rm supp}(f_j)$.  Clearly, then, there are either at least $\kappa \delta_2 N$ pairs
$(b,c) \in B \times C$ for which $T(b,c) = x$ and $b \not \in {\rm supp}(f_j)$; or, 
there are at least $\kappa \delta_2 N$ pairs $(b,c) \in B \times C$ for which
$T(b,c) = x$ and $c \not \in {\rm supp}(f_j)$.  Suppose that the former holds;
then, since for fixed $b$ and $x$ there can be at most $\kappa$ choices 
for $c \in G$ with $T(b,c) = x$, it follows that there are at least $\delta_2 N$ elements $b \in B$
that do not belong to ${\rm supp}(f_j)$.  And if the latter holds, then there are at least $\delta_2 N$
elements $c \in C$ that do not belong to ${\rm supp}(f_j)$.  
Clearly, then, the support of $f_{j+1}$ is larger than the support of $f_j$ by at least $\delta_2 N$ elements.
Iterating this, and using the fact that 
$$
|A|\ =\ {\rm supp}(1_A)\ =\ \theta N\ \geq\ \delta_2 N,
$$ 
we deduce that the support of $f_j$ has size at least $j \delta_2 N$, 
which implies that the procedure must terminate with a
function $f_J$ where $J \leq \delta_2^{-1}$.  We then just let $f := J^{-1} f_J$.  And now, 
since 
$$
f_J*_T f_J(x)\ =\ J^2 f*_T f(x)\ \leq\ \delta_2^{-2} f*_T f(x)\ {\rm  for\ every\ }x \in G,
$$  
and since at this last iteration that produced $f_J$ we stopped at step 2 in the above algorithm, 
it follows that for every pair of sets $B,C \subseteq G$ such that $\|1_A - 1_B\|, \|1_A - 1_C\| \leq \delta_1$,
$$
1_B*_T 1_C(x)\ \leq\ \delta_2^{-2} f*_T f(x) + 2\kappa \delta_2 N\ {\rm for\ every\ } x \in G.
$$

It remains to show that $\| 1_A - f\| \leq \delta_1$:  we first note that $f = J^{-1}(1_{B_1} + \cdots + 1_{B_J})$,
where $B_1,...,B_J$ are the sets $B$ or $C$ arising at each iteration of step 3 in the above
algorithm, and satisfy $\|1_{B_i} - 1_A\| \leq \delta_1$.  Then, writing
$$
f - 1_A\ =\ J^{-1}(( 1_{B_1} - 1_A) + \cdots + (1_{B_J} - 1_A))
$$ 
the triangle inequality immediately gives us that $\|f - 1_A\| \leq \delta_1$,
thereby completing the proof of the lemma.

\section{Proof of main theorems} \label{proofs_section}

\subsection{Proof of Theorem \ref{A_convolve_theorem}}

Let $S_0$ be a random subset of $G$ of size
$$
K\ :=\ \lfloor c \eps^{-2} \delta^{-6} \theta^{-10} 
(\log N - \log(\delta \eps \theta))\rfloor,
$$
where $c > 0$ is some constant to be determined later.

We will show that with positive probability the set $S_0$
satisfies the conclusion of the theorem (with $S_0$ playing the role of $S$) 
for every set $A$ satisfying $|A| = \theta_0 N$
and for every possible set density $\theta_0 \geq \theta$.  And therefore, there
exists a set $S$ with the properties claimed by the theorem.

So, let $\theta_0 \geq  \theta$ be an arbitrary set density, and let
$\delta_1 = \delta^3 \theta_0^{5.5} \sqrt{\eps}/128$ and $k = \lfloor 4\delta_1^{-2} \theta_0 \rfloor +1$.  
Associate to a subset $C \subseteq G$, $|C| = \theta_0 N$, a vector 
$$
v_C\ :=\ (\chi_1,...,\chi_k,\Re \hat 1_C(\chi_1), \Im \hat 1_C(\chi_1), ..., \Re \hat 1_C(\chi_k), \Im \hat 1_C(\chi_k)).
$$ 
where $\chi_1,...,\chi_k$ are the places corresponding to the $k$ largest Fourier coefficients of $h = 1_C$,
as described in Lemma \ref{Fourier_upper_bound}.  From the conclusion of that same lemma we deduce that
\begin{equation} \label{hatbound}
|\hat 1_C(\chi_k)|\ \leq\ \delta_1 N/2.
\end{equation}

Next, round the last $2k$ coordinates of $v_C$ to the nearest multiple of
$\delta_1 N/2$.  Let $w_C$ denote the new vector that results.  It is obvious that as we vary over 
subsets $C \subseteq G$ satisfying $|C| = \theta_0 N$,  the number of possibilities for $w_C$ 
is bounded (crudely) from above by $N^k (\delta_1/4)^{-2k}$.

Furthermore, if two sets $B$ and $C$, $|B|=|C|=\theta_0 N$ share the same vector $w_B = w_C$ 
then the last $2k$ coordinates of $v_B$ and $v_C$ come within 
$\delta_1 N/2$ of one another; and, in light of (\ref{hatbound}), this implies that
$$
\max_{\chi \in \hat G} |\hat 1_B(\chi) - \hat 1_C(\chi)|\ \leq\ \delta_1 N.
$$

The possibilities for the vectors $w_C$ give us a way of placing sets $C \subseteq G$
with $|C| = \theta_0 N$ into equivalence classes.  And now suppose we have one of
these equivalence classes containing a set $C$ such that $|C+C| \leq (1-\eps)N$.  

Applying Corollary \ref{main_lemma}, with the
role of $A$ in that corollary played by our set $C$, the role of $C$ played by our set $B$,
and the role of $\theta$ played by $\theta_0$,
we find that for $\delta_2 = \delta \theta_0^2/4$ there exists a function 
$f : G \to [0,1]$ such that $f - 1_C$ is $\delta_1$-uniform
and such that for every set $B$ where $1_B - 1_C$ is $\delta_1$-uniform we have for every
$x \in G$,
$$
1_B*1_B(x)\ \leq\ \delta_2^{-2} f*f(x) + 2 \delta_2 N.
$$
Let now $U$ denote the set of $x \in G$ where $f*f(x) \leq \delta^3 \theta_0^6 N/32$.  Then, for
every $x \in U$ we have 
$$
1_B*1_B(x)\ <\ \delta \theta_0^2 N.
$$

Note that the fact this holds for all sets $B$ such that $1_B - 1_C$ is $\delta_1$-uniform 
implies that it holds for all sets $B$ in the same equivalence class (described earlier in the proof) 
as $C$.

We next show that this set $U$ is ``large'', by showing that it contains many elements of 
$G \setminus (C+C)$:   since $f - 1_C$ is $\delta_1$-uniform, from
Lemma \ref{convolutions_uniformity_lemma}, using $g = 1_C$ and $h = f$, we have that
$$
\sum_{x \in \Z_N} |f*f(x) - 1_C*1_C(x)|^2\ \leq\ \delta_1^2 (4\theta_0 + \delta_1) N^3. 
$$
It follows that if we let $T$ denote the set of all $x \in G$ where $1_C*1_C(x) = 0$
and where $f*f(x) > \delta^3 \theta_0^6 N/32$, then
$$
|T| \delta^6 \theta_0^{12} N^2/1024\ \leq\ \delta_1^2 (4\theta_0 + \delta_1) N^3;
$$
that is,
$$
|T|\ \leq\ 1024 \delta_1^2(4\theta_0 + \delta_1) \delta^{-6} \theta_0^{-12} N.
$$
Since $|C+C| \leq (1-\eps)N$  it follows that there are at
least $(\eps - 1024 \delta_1^2 (4\theta_0 + \delta_1) \delta^{-6} \theta_0^{-12})N \geq \eps N/2$
places $x \in U$.

And now, because there are at most $N^k (\delta_1/4)^{-2k}$ possible vectors $w_C$,
and therefore at most that many equivalence classes, 
it follows that we have a collection of at most $N^k (\delta_1/4)^{-2k}$ sets $U$ such that for 
every set $B \subseteq G$ satisfying $|B+B| \leq (1-\eps)|B|$, the collection contains 
a set $U = U_B$, $|U| \geq \eps N/2$, satisfying 
$$
1_B*1_B(x)\ <\ \delta \theta_0^2 N,\ {\rm for\ every\ }x \in U.
$$
We claim that with positive probability our random set $S_0$ (defined at the beginning of
the proof) will have non-empty intersection with all these 
$< N^k (\delta_1/4)^{-2k}$ sets $U$.  
To see this, first note that for $c > 0$ sufficiently large in the definition
of $K$ above, the probability that $S_0$ fails to intersect any 
particular set $U_B$ is at most 
$$
{{(1-\eps/2)N \choose K} \over {N \choose K}}\ \leq\ (1-\eps/2)^K\ \leq\ \exp(-K \eps/2)\ <\ N^{-k-1} (\delta_1/4)^{2k}.
$$
So, by the union bound, the probability that $S_0$ intersects {\it every} one of our sets $U_B$ is positive, even
as we vary over all the $\leq N$ choices for set density $\theta_0 \geq \theta$.

It follows that there exists a set $S$ of size $K$ that intersects all the sets $U_B$; and therefore if
$A \subseteq G$, $|A| = \theta_0 N \geq \theta N$, and if $1_A*1_A(x) > \delta \theta_0^2 N$ for every $x \in S$, 
$A$ could not belong to an equivalence class containing a set $C$ such that $|C+C| \leq (1-\eps)N$.
In particular, this would mean that $|A+A| > (1-\eps)N$, thereby completing the proof of the theorem.

\subsection{Proof of Theorem \ref{A_sum_limits}}

The proof of this theorem is not much more than the Dirichlet Box Principle.  Before we state it, we
introduce the notation $\| x\|$ to denote the least residue mod $N$ in absolute value that is 
congruent to $x$ mod $N$. 

\begin{lemma}  Suppose that $x_1,...,x_k \in \Z_N$.  Then, there exists an integer 
$n \not \equiv 0 \pmod{N}$ satisfying 
$$
\| n x_1 \|,\ \|n x_2 \|,\ ...,\ \| n x_k\|\ \leq\ N^{1-1/(k+1)}.
$$
\end{lemma}

\noindent {\bf Proof of the lemma.}  The proof of this lemma is standard:  we consider the set of 
vectors 
$$
\{ (j x_1,\ ...,\ j x_k) \pmod N\ :\ 0 \leq j \leq N^{1-1/(k+1)}\}\ \subseteq\ (\R / N \Z)^k.
$$
Around each point draw a $k$-dimensional box (so, the point 
$(jx_1, ..., j x_k) \pmod{N}$ is the center point of the box) having edge
length $N^{1-1/(k+1)}$.  The total volume consumed by all the boxes exceeds
$$
N^{k-k/(k+1)} N^{1-1/(k+1)}\ =\ N^k.
$$  
So, at least two of those boxes must have a point in common;
say these boxes correspond to $j = a$ and $j = b$, where $a < b$.  It follows
that for each $i=1,...,k$ the numbers $a x_i$ and $b x_i$ are at most $N^{1-1/(k+1)}$
apart when considered mod $N$.  Letting $n = b-a$ it is easy to see that this then 
implies the lemma. \hfill $\blacksquare$
\bigskip

Given the set $S$, set $k = |S|$ and $\{x_1,...,x_k\} = S$, and then let $n$ be as in the lemma.
We then let
$$
A\ :=\ n^{-1}*B\ \subseteq\ \Z_N,\ {\rm where\ } B\ :=\ \{t\ :\ -N/6 < t < N/6\}\ \subseteq\ \Z_N,
$$
where the notation $\lambda * B$ represents the set that results when 
we dilate the elements of $B$ by $\lambda$.

It follows that for each $x \in S$ we have 
$$
1_A*1_A(x)\ =\ 1_{n^{-1}*B}*1_{n^{-1}*B}(x)\ =\ 1_B*1_B(n x).
$$
We have that if $k < (\log N)/2$ then $\| nx\| \lesssim N/e^2$ for $x\in S$, and then 
$$
1_A*1_A(x)\ >\ 1_B*1_B(\lfloor N/e^2 \rfloor+1)\ \gtrsim\ N(1/3 - 1/e^2)\ >\ N/6.
$$
This completes the proof of Theorem \ref{A_sum_limits}.

\subsection{Proof of Theorem \ref{A_sum_pseudo}}

To prove the theorem we will show that if $A \subseteq G$ is any
set satisfying $|A+A| < (1-\eps)|A|$, where $|A| = \theta N$, then the level-set 
$\{x \in G\ :\ 1_A*1_A(x) \leq \delta \theta^2 N\}$ must have non-trivial intersection with the set $S$.
And, from Theorem \ref{longaps} we furthermore have that to prove this conclusion 
it suffices to show that $S$ intersects every translate of every Bohr neighborhood of 
radius $\delta^3 \theta^{4.5} \sqrt{\eps}/128$ and dimension at most 
$2^{16}\delta^{-6} \theta^{-10} \eps^{-1} + 1$.

Letting $\B(\Lambda, \rho)$ be such a Bohr neighborhood, where $\rho = \delta^3 \theta^{4.5} \sqrt{\eps}/128$,
we begin by letting $\B' = \B(\Lambda, \rho/2)$ and setting $g = 1_{\B'}*1_{\B'}$.  Since 
${\rm supp}(g) \subseteq \B(\Lambda, \rho)$, to prove our theorem it suffices to show that for 
every translate $t \in G$ we have
\begin{equation} \label{sumgs}
\sum_{x \in G} g(x+t) 1_S(x)\ >\ 0.
\end{equation}
In terms of Fourier coefficients this inequality is simply
$$
\sum_{\chi \in \hat G} \chi(t) \hat 1_{\B'}(\chi)^2 \hat 1_S(\chi^{-1})\ >\ 0;
$$
and to prove this it suffices to show that
$$
\max_{\chi \neq \chi_0} |\hat 1_S(\chi)| \sum_{\chi \neq \chi_0} |\hat 1_{\B'}(\chi)|^2 \ <\ |\B'|^2 |S|.
$$
From Parseval's identity we have that this holds provided
$$
\max_{\chi \neq \chi_0} |\hat 1_S(\chi)|/|S|\ <\ |\B'|/N.
$$

To finish our proof we apply Lemma \ref{bohr_lemma}
using $r = \rho/2$, and deduce that we have that (\ref{sumgs}) holds provided 
$$
\max_{\chi \neq \chi_0} |\hat 1_S(\chi)|/|S|\ <\ 
(\delta^3 \theta^{4.5} \sqrt{\eps}/512 \pi)^{2^{16} \delta^{-6} \theta^{-10} \eps^{-1}  + 1},
$$
which is one of our assumptions.  The theorem now follows.

\subsection{Proof of Theorem \ref{longaps}}

Arrange the Fourier coefficients of $h = 1_A$ from largest to smallest in magnitude
as in Lemma \ref{Fourier_upper_bound}.

Let $\delta_1 =  \delta^3 \theta^{5.5} \sqrt{\eps}/128$ and let
$k = \lfloor 4\delta_1^{-2} \theta \rfloor +1$.  Letting $\chi_1,...,\chi_N$ be
as in (\ref{respect_me}), from Lemma \ref{Fourier_upper_bound}
we will have that $|\hat 1_A(\chi_k)|\ \leq\ \delta_1 N/2$.

Next, we let $\B$ denote the Bohr neighborhood $\B(\chi_1,...,\chi_k; \delta_1/\theta)$. 
Then, for each $t \in \B$ define the set $A_t := A+t$.  Note that $1_{A+t}(x) = 1_A(x-t)$ and that
$\hat 1_{A_t}(\chi) = \chi(t) \hat 1_A(\chi)$.

We have that for $i=1,2,...,k$,
$$
|\hat 1_{A_t}(\chi_i) - \hat 1_A(\chi_i)|\ \leq\ |1 - \chi_i(t)|\cdot |\hat 1_A(\chi_i)|\ \leq\ \delta_1 N.
$$
And for $k+1 \leq i \leq N$ we have
$$
|\hat 1_{A_t}(\chi_i) - \hat 1_A(\chi_i)|\ \leq\ |1 - \chi_i(t)| \cdot |\hat 1_A(\chi_i)|\ \leq\ 2 |\hat 1_A(\chi_i)|\ \leq\ \delta_1 N.
$$
In general, then, we have that for all $\chi \in \hat G$,
$$
|\hat 1_{A_t}(\chi) - \hat 1_A(\chi)|\ \leq\ \delta_1 N,
$$
which implies that $1_{A_t} - 1_A$ is $\delta_1$-uniform.

We now apply Corollary \ref{main_lemma} using $\delta_2 = \delta \theta^2/4$, and
deduce the existence of a function $f : G \to [0,1]$ having the properties indicated by the Corollary.
From the fact that $f - 1_A$ is $\delta_1$-uniform we deduce from Lemma \ref{convolutions_uniformity_lemma},
using $g = 1_A$ and $h= f$, that
\begin{equation} \label{ffaa}
\sum_{x \in G} |1_A*1_A(x) - f*f(x)|^2\ =\ \delta_1^2 (4\theta + \delta_1) N^3.
\end{equation}
Letting $T$ denote the set of all $x \in G$ such that $1_A*1_A(x) < \delta^3 \theta^6 N/128$,
since we are given that $|T| \geq \eps N$ we deduce from (\ref{ffaa}) that if for all such
$x$ we also had that $f*f(x) \geq \delta^3 \theta^6 N/32$, then
$$
9 \delta^6 \theta^{12} \eps N^3/2^{14}\ \leq\ |T| (3\delta^3 \theta^6 N/128)^2\ \leq\ \delta_1^2 (4\theta + \delta_1) N^3,
$$
which is impossible.  So, there exists $x \in G$ such that $f*f(x) < \delta^3 \theta^6 N/32$.
From Corollary \ref{main_lemma}, and
the fact that $1_A - 1_{A_t}$ is $\delta_1$-uniform (and that $1_A - 1_A = 0$ is also
$\delta_1$-uniform), we will have for this value $x$ that
$$
1_A*1_A(x-t)\ =\ 1_A*1_{A_t}(x)\ \leq\ \delta_2^{-2} f*f(x) + 2 \delta_2 N\ <\ \delta \theta^2 N.
$$
It follows that the level-set $\{y\ :\ 1_A*1_A(y) < \delta \theta^2 N\}$ contains the set $x - \B$,
which completes the proof of the theorem.

\section{Acknowledgments} \label{ack}

I would like to thank Olof Sisask and Thomas Bloom for their numerous comments and suggestions.
I would  also like to thank the referee for the comment about using Green's regularity lemma and
Bohr neighborhoods.


\begin{thebibliography}{999}

\bibitem{ajtai} M. Ajtai, H. Iwaniec, J. Koml\'os, J. Pintz, and E. Szemer\'edi, {\it Construction of a thin set with
small Fourier coefficients}, Bull. London Math. Soc. {\bf 22} (1990), 583-590.

\bibitem{bateman} M. Bateman and N. H. Katz, {\it New bounds on cap sets}, preprint arXiv:1101.5851

\bibitem{behrend} F. A. Behrend, {\it On sets of integers which contain no three terms in arithmetic progression},
Proc. Nat. Acad. Sci. {\bf 23} (1946), 331-332.

\bibitem{bloom} T. Bloom, {\it Translation invariant equations and the method of Sanders}, preprint arXiv:1107.1110

\bibitem{bogolyubov} N. Bogolio\` uboff, {\it Sur quelques propi\'et\'es arithm\'etiques des presque-p\'eriodes},
Ann. Chaire Phys. Math. Kiev {\bf 4} (1939), 185-205.

\bibitem{bourgain} J. Bourgain, {\it On triples in arithmetic progression}, Geom. Funct. Anal. {\bf 9} (1999), 968-984.

\bibitem{croot} E. Croot, I. Laba, and O. Sisask, {\it Arithmetic progressions in sumsets and $L^p$-almost-periodicity},
preprint arXiv:1103.6000

\bibitem{croot2} E. Croot and O. Sisask, {\it A probabilistic technique for finding almost-periods of convolutions},
Geom. Funct. Anal. {\bf 20} (2010), no. 6, 1367-1396.

\bibitem{elkin} M. Elkin, {\it An improved construction of progression-free sets} Israel Jour. of Math. {\bf 184} (2011), 93-128.

\bibitem{gowers} W. T. Gowers, {\it A new proof of Szemer\'edi's Theorem for arithmetic progressions of length four}, 
Geom. Funct. Anal. {\bf 8} (1998), 529-551.

\bibitem{green} B. Green, {\it Arithmetic progressions in sumsets}, Geom. Funct. Anal. {\bf 12} (2002), no. 3, 584-597.

\bibitem{green_regular} ---------------, {\it A Szemer\' edi-type regularity lemma in abelian groups, with applications},
Geom. Funct. Anal {\bf 15} (2005), 340-376.

\bibitem{wolf} B. Green and J. Wolf, {\it A note on Elkin's improvement of Behrend's construction}, Additive Number 
Theory: Festschrift in Honor of the Sixtieth Birthday of Melvyn B. Nathanson, 2010, p. 141-144. 

\bibitem{katz} N. H. Katz and P. Koester, {\it On additive doubling and energy}, SIAM J. Discrete Math. {\bf 24} (2010), 1684-1693.

\bibitem{roth} K. F. Roth, {\it On certain sets of integers}, J. London Math. Soc. {\bf 28} (1953), 104-109.

\bibitem{ruzsa} I. Z. Ruzsa, {\it Arithmetic progressions in sumsets}, Acta Arith. {\bf 60} (1991), no. 2, 191-202.

\bibitem{sanders} T. Sanders, {\it On the Bogolyubov-Ruzsa Lemma}, preprint arXiv:1011.0107

\bibitem{sanders2} ---------------, {\it On Roth's Theorem on progressions}, to appear in Ann. of Math.

\bibitem{schoen} T. Schoen and I. Shkredov, {\it Roth's Theorem in many variables}, preprint arXiv:1106.1601

\bibitem{taovu} T. Tao and V. Vu, {\it Additive Combinatorics}, 2006 Cambridge Univ. Press.

\end{thebibliography}
\end{document}